\documentclass[journal]{IEEEtran}
\usepackage{graphicx}
\usepackage{cite}
\usepackage{picinpar}
\usepackage{amsmath}
\usepackage{mathtools}
\DeclarePairedDelimiter\norm{\lVert}{\rVert}
\usepackage{url}
\usepackage{flushend}
\usepackage[latin1]{inputenc}
\usepackage{colortbl}
\usepackage{soul}
\usepackage{multirow}
\usepackage{pifont}
\usepackage{color}
\usepackage{alltt}
\usepackage[hidelinks]{hyperref}
\usepackage{enumerate}
\usepackage{url}
\usepackage{epstopdf}
\usepackage{pbox}
\usepackage{fancyhdr}
\usepackage{lineno} 
\usepackage{bm}

\usepackage{amssymb}
\usepackage{stfloats}
\usepackage{algorithmic}
\usepackage[ruled,vlined]{algorithm2e}

\usepackage{subfigure}
\usepackage{caption}

\usepackage[table]{xcolor}
\usepackage{comment}
\newcommand{\rv}[1]{\textcolor{black}{#1}} 
\begin{document}

\title{Network-Aware Asynchronous Distributed ADMM Algorithm for Peer-to-Peer Energy Trading}

\author{
Zeyu~Yang and Hao~Wang

\thanks{Z. Yang is with the Department of Computer Science, Rice University, Houston, TX 77005, USA (e-mail: zeyu.yang@rice.edu).}
\thanks{H. Wang is with the Department of Data Science and Artificial Intelligence, Faculty of Information Technology and the Monash Energy Institute, Monash University, Melbourne, VIC 3800, Australia (e-mail: hao.wang2@monash.edu).}
}

\maketitle

\begin{abstract}
The increasing uptake of distributed energy resources (DERs) in smart home prosumers calls for distributed energy management strategies, and the advances in information and communications technology enable peer-to-peer (P2P) energy trading and transactive energy management. 
Many works attempted to solve the transactive energy management problem using distributed optimization to preserve the privacy of DERs' operations. 
But such distributed optimization requires information exchange among prosumers, often via synchronous communications, which can be unrealistic in practice.
This paper addresses a transactive energy trading problem for multiple smart home prosumers with rooftop solar, battery storage, and controllable load, such as heating, ventilation, and air-conditioning (HVAC) units, considering practical communication conditions. We formulate a network-aware energy trading optimization problem, in which a local network operator manages the network constraints supporting bidirectional energy flows. We develop an asynchronous distributed alternating direction method of multipliers (ADMM) algorithm to solve the problem under asynchronous communications, allowing communication delay and indicating a higher potential for real-world applications. We validate our design by simulations using real-world data. The results demonstrate the convergence of our developed asynchronous distributed ADMM algorithm and show that energy trading reduces the energy cost for smart home prosumers.
\end{abstract}


\markboth{IEEE}%
{}

\definecolor{limegreen}{rgb}{0.2, 0.8, 0.2}
\definecolor{forestgreen}{rgb}{0.13, 0.55, 0.13}
\definecolor{greenhtml}{rgb}{0.0, 0.5, 0.0}
\SetAlFnt{\small}
\SetAlCapFnt{\small}

\section{Introduction}\label{sec:intro}
The recent global decarbonization trend promotes energy efficiency for home energy use, such as heating, ventilation, and air-conditioning (HVAC), and the adoption of sustainable energy devices, such as rooftop solar photovoltaic (PV) and battery storage. 
Energy management for flexible energy sources is essential to improve energy efficiency and reduce energy costs. For example, HVAC units consume a large amount of energy in residential households and buildings, but their demand is flexible. The studies in \cite{yoon2018optimal,mansy2020optimal} proposed adjusting HVAC energy consumption to perform demand response in a time-varying environment, such as dynamic electricity prices and intermittent renewables. Wu and Skye in \cite{wu2018net} studied the configuration of HVAC and solar systems sized to reach the net-zero transition for residential buildings. 
Recent studies in \cite{li2018real,arun2017intelligent,paul2020real} proposed optimization strategies to minimize costs by energy scheduling for smart homes equipped with battery storage, rooftop solar, HVAC, and other smart home appliances.

In addition to the energy scheduling of smart homes, cooperative energy management and energy trading have been emerging as a new paradigm enabled by two-way communications, such as Internet of things (IoT) technologies. 
Smart homes and their distributed energy resources (DERs) can communicate with each other to coordinate energy scheduling and trading by leveraging their diverse generation and demand profiles. For example, the optimal energy management for a cluster of residential prosumers was studied in \cite{luna2016cooperative} considering cooperative operations among prosumers, such as energy sharing. 
Alam et al. in \cite{alam2019peer} studied an energy cost optimization problem by coordinating peer-to-peer (P2P) energy trading among smart homes in a demand-side management system. Energy trading further improves the smart grid operation by utilizing solar generation and leveraging the flexibility of energy consumption. However, these studies adopted a centralized solution requiring the information of all smart homes for decision making, thus causing concerns about privacy and feasibility due to the diverse ownership of DERs in smart home prosumers.

Distributed energy management and energy trading are becoming promising solutions to preserve privacy while facilitating the effective coordination of DERs. For example, a multi-leader multi-follower noncooperative Stackelberg game was formulated in \cite{mondal2017demands} to enable distributed energy management among connected customers in microgrids. 
A distributed energy management method was developed in \cite{yang2021distributed} for P2P energy trading among smart homes with HVAC systems and rechargeable batteries to reduce the energy cost and preserve privacy. 
A study in \cite{yang2021privacy} developed a transactive energy management algorithm using the alternating direction method of multipliers (ADMM) for smart homes to trade energy. An ADMM distributed algorithm was proposed in \cite{kou2020scalable} for coordinating residential devices at scale to enable demand response. \rv{We include a comprehensive literature review and analysis on transactive energy and energy trading in Section~\ref{sec:related}, covering centralized, decentralized, and distributed solution methods.}
Most of the distributed energy management and energy trading methods, including ADMM-based methods, assumed synchronous communications for DERs to solve the energy management problem. However, due to the heterogeneous communication capabilities of DERs and latency in transmitting energy trading information, it is more realistic and practical to develop distributed energy management strategies for DERs under asynchronous communications. 

To bridge the research gap, this paper develops a cooperative energy management system for energy trading among smart home prosumers under asynchronous communications. 
We formulate energy management and energy trading among prosumers as a total cost minimization problem. A local network operator manages the distribution feeder to ensure the two-way energy flows satisfy the network constraints.
We develop an asynchronous distributed ADMM algorithm to solve the optimal energy management and energy trading to optimize smart homes' thermal comfort and energy cost. Unlike traditional distributed algorithms, our developed asynchronous distributed ADMM algorithm does not require synchronous decision making and information exchange but converges under asynchronous communications.
We validate our design through simulations using real-world data, and the results demonstrate the convergence of our developed asynchronous distributed ADMM algorithm.
The numerical results also show that even with homogeneous solar generations in a community region, energy trading can achieve a cost reduction of up to $34\%$ for prosumers by exploiting the flexibility in HVAC energy consumption and battery storage.

\section{\rv{Related Works}}\label{sec:related}
We review related works on energy management and P2P energy trading from the perspective of information exchange and computing in this Section.

Effective energy management is at the heart of the integration and orchestration of DERs. Among various controllable loads on the user side, HVAC consumes a large amount of energy, and thus the energy optimization of HVAC systems has attracted a lot of attention. For example, Yu et~al. \cite{yu2017distributed} developed a real-time control method to minimize the cost of an HVAC system using Lyapunov optimization techniques. 
The energy efficiency of HVAC systems in multi-zone buildings was studied in \cite{song2020energy} to enable personalized control strategies.
Pairing HVAC systems with PV systems was studied in \cite{wu2018net} to facilitate the net-zero transition for buildings. Energy storage \cite{zhang2021review} also plays a vital role in energy systems for integrating DERs and improving the system performance.
Instead of optimizing a single or aggregated HVAC system, recent studies focused on leveraging energy trading as a promising smart grid technology to coordinate HVAC systems. For example, a close-to-optimal energy trading method was studied in \cite{ALAM20191434} to minimize the electricity cost of smart homes, but the method is centralized. To solve such a centralized problem, the controller requires all the information about HVAC systems, leading to serious privacy concerns. According to the study in \cite{wang2020activity}, the HVAC energy consumption contains rich information of occupants' activities, and thus privacy-preserving is essential in energy management solutions.
Recent studies in \cite{yang2021distributed} developed a privacy-preserving distributed algorithm to optimize the HVAC energy consumption of cooperative smart homes via energy trading. There are some challenges to HVAC energy management, e.g., the baseline estimation in HVAC demand response and the optimization of pre-cooling based on occupancy detection, which needs further study.

Instead of focusing on the HVAC system, more recent research attempted to optimize the DER energy management. For example, a centralized control was proposed in \cite{caldognetto2014improving} to reduce the power loss of DERs and enhance the network efficiency. A centralized solution was studied in \cite{luna2016cooperative} to coordinate DERs for sharing energy among smart homes. Game-theoretic approaches have also been used to model and solve the optimal interactions between the network and DERs. A Stackelberg game was formulated in \cite{liu2017energy} to optimize the operation of microgrids with PV prosumers and effective energy sharing management. Similarly, cooperative game models, such as the coalition game, were used in \cite{mei2019coalitional} to provide incentives to coalitional operations of microgrids for power exchange. 
However, the above centralized solution methods require complete information of all DERs, causing privacy concerns. Therefore, recent studies aimed to develop privacy-preserving algorithms for energy management and energy trading of DERs.

From the perspective of information exchange and computing, related works on energy trading can be divided into centralized, decentralized, and distributed solution methods. As discussed above, centralized methods cannot preserve privacy, but decentralized and distributed methods promise privacy preservation. Enabled by information and communications technology (ICT), such as IoT, DERs can exchange information with each other. Hence, with proper designs of necessary information exchange protocols, DERs can cooperatively solve energy management and trading problems in decentralized and distributed manners. For example, decentralized solution methods have been proposed in \cite{yang2019fully,wang2021consensus}. A fully decentralized energy market at the distribution level was developed in \cite{yang2019fully} without cooperation among market participants. Wang et al. in \cite{wang2021consensus} developed a fully decentralized transactive energy management system using a consensus-based algorithm and designed a virtual pool for DER prosumers to trade energy without revealing private information. 
However, one of the challenges to decentralized methods is that prosumers or DERs do not have global information, particularly about the energy network. Thus the role of the network operator or aggregator is inevitable in enforcing the network constraints.

Different from decentralized methods discussed above, distributed methods leverage the role of the network operator to satisfy the network constraints for energy trading while preserving the privacy for DERs. Primal-dual decomposition and ADMM are widely used techniques to design distributed algorithms. For example, a distributed optimization algorithm was designed in \cite{rokni2018optimum} for energy management of supplies and demands in a microgrid using ADMM. A cooperative bargaining model was developed to incentivize energy trading among interconnected microgrids and solved by a privacy-preserving ADMM distributed algorithm \cite{wang2016incentivizing}. A dual-decomposition-based distributed algorithm was proposed in \cite{qu2020distributed} for a micro-integrated energy system, in which electricity is delivered to customers by energy hubs. ADMM was used in \cite{mu2020decentralized} to solve the economic dispatch of micro-integrated energy systems for privacy preservation. However, these primal-dual decomposition and ADMM algorithms usually require synchronous communications among DERs. In other words, the computation in these algorithms is conducted in a synchronous manner, waiting for information exchange to be completed by all DERs. Due to the heterogeneous nature of DERs, the communication capability, bandwidth, and latency can be very different, resulting in asynchronous information exchange in the distributed algorithm. But existing works lack such a consideration to accommodate asynchronous information exchange. This motivates our work to consider this practical need and develop an asynchronous ADMM distributed algorithm for energy trading under asynchronous communications.

\section{System Model And Problem Formulation}\label{sec:model}
\begin{figure}[t]
    \centering
    \includegraphics[width=1.0\linewidth]{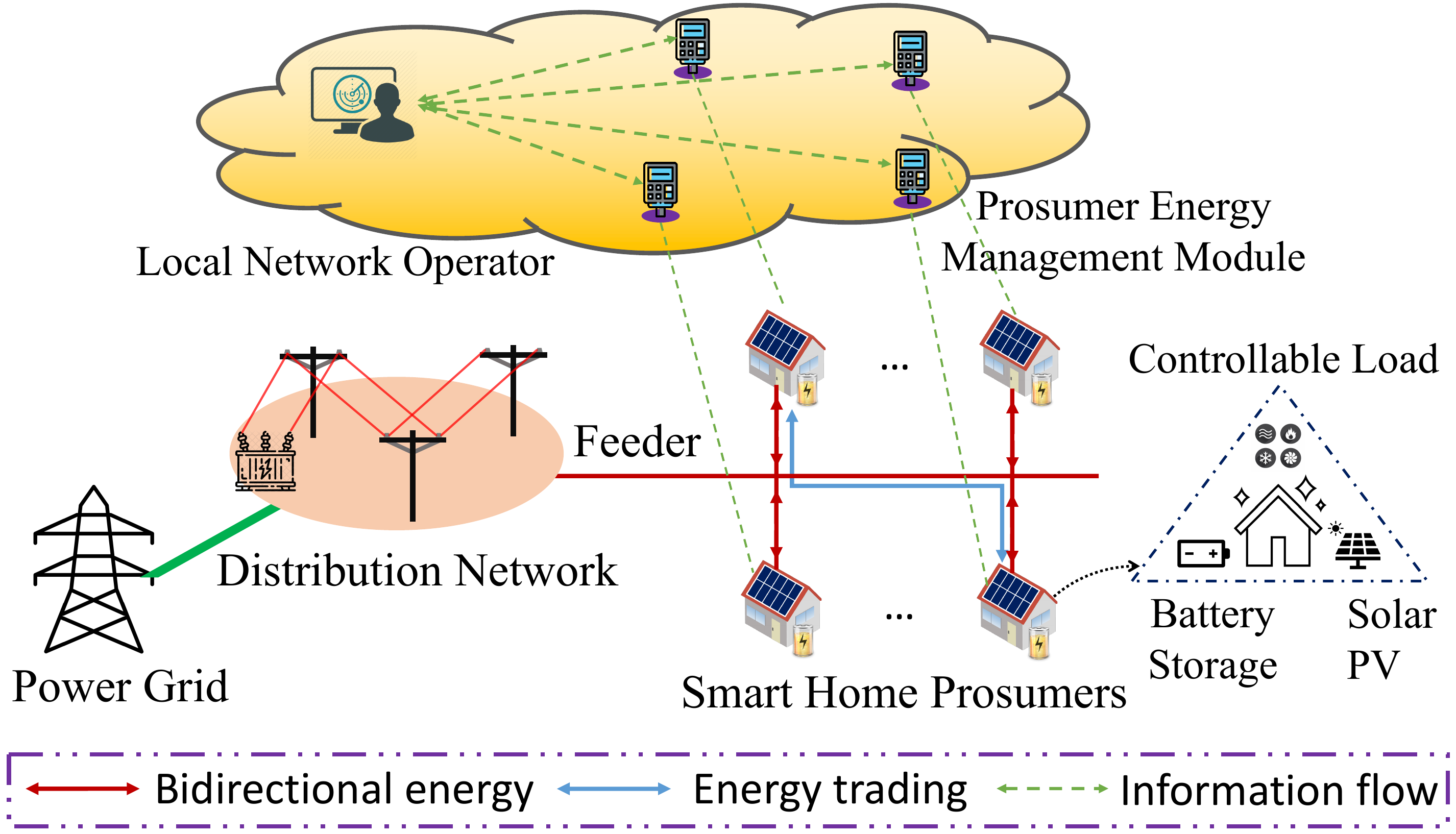}
    \caption{\rv{The system diagram for energy management and energy trading of smart home prosumers in an energy network with information exchange enabled by two-way communications.}}
    \label{f:problem}
\end{figure}

We consider an energy management and energy trading system for smart home prosumers denoted as $\mathcal{N} {=} \{1,..., N\}$ in a community, as shown in Fig. \ref{f:problem}. Smart home prosumers aim to optimize their energy management of PV, battery, and controllable load, such as HVAC, and energy trading with other prosumers to minimize the energy cost over an operational horizon $\mathcal{H} {=} \{1,..., H \}$. Note that smart home prosumers are connected to the low-voltage power distribution grid to trade energy and can exchange information using two-way communications. \rv{The feeder can support bidirectional energy flow, e.g., prosumers can feed in excessive solar PV energy to the grid and trade energy with peer prosumers.} Through energy trading, prosumers can leverage the flexibility and diversity of HVAC energy consumption and battery scheduling for mutual benefits.

\subsection{Model for Energy Management and Energy Trading}\label{sec:smarthome}
We model the energy supply and consumption in smart home prosumers as follows. Since we focus on the day-ahead energy management and energy trading, we will define all decision variables and parameters in vector forms for a length of a day with $H$ time slots.

\subsubsection{Energy Supply Model}\label{sec:supply}
Smart home prosumers have rooftop solar PV, which generates solar energy $\overline{\bm{P}}_i^{\text{S}} = \{\overline{P}_i^{\text{S}}[t],~\forall t \in \mathcal{H}\}$, where $\overline{P}_i^{\text{S}}[t]$ denotes the solar energy generation at smart home prosumer $i$ in time slot $t$. Smart home prosumers will use solar energy as much as possible to meet their demand and trade with others, and the used solar energy is denoted as $\bm{p}_i^{\text{S}} = \{p_i^{\text{S}}[t],~\forall t \in \mathcal{H}\}$. \rv{Prosumers can feed in excessive solar energy to the grid, denoted as $\bm{p}_i^{\text{feedin}} = \{p_i^{\text{feedin}}[t],~\forall t \in \mathcal{H}\}$, and earn revenue from selling it. But due to the local network constraints, prosumers may not feed in all excessive solar energy, and the remaining solar energy will be curtailed.} The solar energy use $\bm{p}_i^{\text{S}}$ and the feed-in solar energy $\bm{p}_i^{\text{feedin}}$ satisfy the following constraint
    \begin{align}
        & \rv{\bm{p}_i^{\text{S}} \geq \bm{0},~\bm{p}_i^{\text{feedin}} \geq \bm{0}}, \label{constraint-solar1} \\
        & \rv{\bm{p}_i^{\text{S}} + \bm{p}_i^{\text{feedin}} \leq \overline{\bm{P}}_i^{\text{S}}}, \label{constraint-solar2}
    \end{align}
where both solar energy use $\bm{p}_i^{\text{S}}$ and the feed-in solar energy $\bm{p}_i^{\text{feedin}}$ are non-negative, and their sum is no greater than the avaiable solar energy $\overline{\bm{P}}_i^{\text{S}}$.

In addition to the local solar energy, prosumers in deficit of local solar energy can also purchase energy from the grid, denoted as $\bm{p}_i^{\text{G}} = \{p_i^{\text{G}}[t],~\forall t \in \mathcal{H}\}$. The grid energy purchase $\bm{p}_i^{\text{G}}$ satisfy the following constraint
    \begin{align}
        \bm{0} \leq  \bm{p}_i^{\text{G}} \leq \overline{\bm{P}}_i^{\text{G}}, \label{constraint-grid}
    \end{align}
where $\overline{\bm{P}}_i^{\text{G}}$ is the vector form of line capacity of smart home prosumer $i$.

\subsubsection{Battery Storage Model}\label{sec:battery}
Smart home prosumers have battery storage units to charge locally generated solar energy and grid energy at low prices and discharge to either meet the demand or sell to peer prosumers as energy trading for mutual benefits.

For prosumer $i$, the charging energy and discharging energy are denoted as $\bm{p}_i^{\text{ch}} = \{p_i^{\text{ch}}[t],~\forall t \in \mathcal{H}\}$ and $\bm{p}_i^{\text{dis}} = \{p_i^{\text{dis}}[t],~\forall t \in \mathcal{H}\}$. The energy stored in the battery storage is denoted as $\bm{p}_i^{\text{B}} = \{p_i^{\text{B}}[t],~\forall t \in \mathcal{H}\}$.
The battery energy dynamics $\Delta p_i^{\text{B}}[t] =  p_i^{\text{B}}[t] - p_i^{\text{B}}[t-1]$ evolves with charging and discharging energy subject to charge and discharge efficiencies, which is shown as 
    \begin{align}
        \Delta \bm{p}_i^{\text{B}} {=} \eta_i^{\text{ch}} \bm{p}_i^{\text{ch}} {-} \frac{\bm{p}_i^{\text{dis}}}{\eta_i^{\text{dis}}}, \label{constraint-battery1} 
    \end{align}
where $\eta_i^{\text{ch}}$ and $\eta_i^{\text{dis}}$ are coefficients for charge and discharge efficiencies.

The battery energy level should be managed within an operational range, often denoted by $\underline{\alpha}_i^{\text{B}} \bm{P}_i^{\text{B}}$ and $\overline{\alpha}_i^{\text{B}} \bm{P}_i^{\text{B}}$, where $\bm{P}_i^{\text{B}}$ is the battery capacity, $\underline{\alpha}_i^{\text{B}}$ and $\overline{\alpha}_i^{\text{B}}$ are the minimum and maximum state of charging. The charging and discharging energy are also restricted by the maximum charging and discharging amount denoted by $\overline{\bm{P}}_i^{\text{ch}}$ and $\overline{\bm{P}}_i^{\text{dis}}$. Therefore, the battery energy level $\bm{p}_i^{\text{B}}$, charging energy $\bm{p}_i^{\text{ch}}$, and discharging energy $\bm{p}_i^{\text{dis}}$ satisfy the following constraints
    \begin{align}
        & \underline{\alpha}_i^{\text{B}} \bm{P}_i^{\text{B}} \leq \bm{p}_i^{\text{B}} \leq \overline{\alpha}_i^{\text{B}} \bm{P}_i^{\text{B}}, \label{constraint-battery2} \\
        & \bm{0} \leq  \bm{p}_i^{\text{ch}} \leq \overline{\bm{P}}_i^{\text{ch}},
         \label{constraint-charge} \\
        & \bm{0} \leq  \bm{p}_i^{\text{dis}} \leq \overline{\bm{P}}_i^{\text{dis}}.
         \label{constraint-discharge}
    \end{align}
        
\subsubsection{Energy Consumption Model}\label{sec:hvac}
Each smart home prosumer has an HVAC unit that consumes energy $\bm{p}_i^{\text{hvac}}= \{p_i^{\text{hvac}}[t],~\forall t \in \mathcal{H}\}$ to control the indoor temperature $\bm{T}_i^{\text{in}}= \{ T_i^{\text{in}}[t],~\forall t \in \mathcal{H}\}$ against outdoor temperature $\bm{T}_i^{\text{out}}= \{ T_i^{\text{out}}[t] ,~\forall t \in \mathcal{H}\}$. The indoor temperature dynamics $\Delta T_i^{\text{in}}[t] =  T_i^{\text{in}}[t] - T_i^{\text{in}}[t-1]$ evolves with the indoor and outdoor temperatures as well as HVAC energy consumption $\bm{p}_i^{\text{hvac}}$ with HVAC efficiency parameters $C_i$, $R_i$ and working modes $\eta_i$. Specifically, the coefficient $C_i$ is the thermal resistance, $R_i$ is the equivalent heat capacity, and $\eta_i$ is the working mode and efficiency, all for prosumer $i$'s HVAC system.
The indoor temperature $\bm{T}_i^{\text{in}}$ and the HVAC energy consumption $\bm{p}_i^{\text{hvac}} \geq \bm{0}$ satisfy the following constraints
    \begin{align}
        & \Delta \bm{T}_i^{\text{in}} = \frac{\bm{T}_i^{\text{in}} - \bm{T}_i^{\text{out}} + \eta_i R_i \bm{p}_i^{\text{hvac}}}{C_i R_i}, \label{constraint-hvac1} \\ %
        & \underline{\bm{T}}_i^{\text{in}} \leq  \bm{T}_i^{\text{in}} \leq  \overline{\bm{T}}_i^{\text{in}}, \label{constraint-hvac2} 
    \end{align}
where Constraint \eqref{constraint-hvac1} shows the dynamics of the indoor temperature changing with the HVAC energy consumption $\bm{p}_i^{\text{hvac}}$. The indoor temperature should be controlled within a comfortable range from $\underline{\bm{T}}_i^{\text{in}}$ to $\overline{\bm{T}}_i^{\text{in}}$, as shown in Constraint \eqref{constraint-hvac2}.

Our work takes the HVAC system as an example for a controllable load. There are also other controllable loads, such as washers and dryers, but these appliances are operated much less frequently compared to HVAC and thus are not considered in our current work. The rest of energy consumption is mostly inflexible, and we denote the inflexible load as $\bm{p}_i^{\text{base}}= \{p_i^{\text{base}}[t],~\forall t \in \mathcal{H}\}$, which is given and not controllable in our model.

\subsubsection{\rv{Energy Trading Model}}\label{sec:trading}
Smart home prosumers can participate in peer-to-peer energy trading behind a common low-voltage transformer in the community. Through ICT, such as IoT communications, \rv{prosumer $i$ can trade energy with other peer prosumers $j \in \mathcal{N} \backslash i$ and decide how much energy to trade, denoted as $\bm{p}_{i}^{\text{trade}} = \{p_{i}^{\text{trade}}[t],~\forall t \in \mathcal{H}\}$. When $p_{i}^{\text{trade}}[t] > 0$, prosumer $i$ purchase energy from peer prosumers; otherwise,  prosumer $i$ sells energy when $p_{i}^{\text{trade}}[t] < 0$.}
Prosumers are sited on the network close to each other, and we assume that the power loss is negligible. The energy trading clearing among all prosumers satisfy the following constraint
    \begin{align}
        \rv{\sum_{i \in \mathcal{N}} \bm{p}_{i}^{\text{trade}} = \bm{0}.} \label{constraint-trading1}
    \end{align}

The energy trading amount for each prosumer is also restricted, and prosumer $i$'s trading capacities are denoted as $\underline{P}_{i}^{\text{trade}} <0$ and $\overline{P}_{i}^{\text{trade}} >0$ for selling and buying, respectively. Hence, the energy trading amount also satisfies the following constraint 
    \begin{align}
        \rv{\bm{\underline{P}}_{i}^{\text{trade}} \leq \bm{p}_{i}^{\text{trade}} \leq \bm{\overline{P}}_{i}^{\text{trade}},} \label{constraint-trading2}
    \end{align}
where $\bm{\underline{P}}_{i}^{\text{trade}}$ and $\bm{\overline{P}}_{i}^{\text{trade}}$ are the vector forms of prosumer $i$'s trading capacities.

Note that, in addition to the explicit energy trading constraints~\eqref{constraint-trading1}~and~\eqref{constraint-trading2}, energy trading also needs to satisfy the local network constraints. We will introduce the network model in Section~\ref{sec:network}. Also, to facilitate the energy trading, the local energy operator will interact with prosumers to derive the optimal decision making through iteration, and we will present the distributed solution method in Section~\ref{sec:solution}.

\subsubsection{Energy Balance in Smart Home Prosumers}\label{sec:balance}
Each smart home prosumer needs to balance the energy supply and demand and thus satisfies the following balancing constraint
    \begin{align}
        \bm{p}_i^{\text{S}} + \bm{p}_i^{\text{G}} + \bm{p}_i^{\text{dis}} + \bm{p}_{i}^{\text{trade}} = \bm{p}_i^{\text{ch}} + \bm{p}_i^{\text{hvac}} + \bm{p}_i^{\text{base}}, \label{constraint-balance} %
    \end{align}
where the left hand side in \eqref{constraint-balance} represents the total energy supply from the solar PV, the grid, battery discharge, and energy trading. The right hand side represents the total energy demand from inflexible load, HVAC load, and battery charge.

\subsection{Model for The Energy Network}\label{sec:network}
For the energy network, we consider a linearized distribution network model, which has been widely used to model distribution network constraints for active power, reactive power, and voltage \cite{zhong2018topology}. 
Each prosumer interact with the grid and peer prosumers through grid energy purchase $\bm{p}_i^{\text{G}}$, feed-in solar energy $\bm{p}_i^{\text{feedin}}$, and energy trading $\bm{p}_i^{\text{trade}}$. Hence, prosumer $i$'s output energy $\bm{p}_i^{\text{output}}= \{p_i^{\text{output}}[t],~\forall t \in \mathcal{H}\}$ is
\begin{align}
    \rv{\bm{p}_i^{\text{output}} = \bm{p}_i^{\text{feedin}} - \bm{p}_i^{\text{G}} - \bm{p}_i^{\text{trade}}.} ~\label{constraint-output}
\end{align}

The linearized power network is modeled as follows. We denote $\bm{p}_i^{\text{injection}} = \{p_i^{\text{injection}}[t],~\forall t \in \mathcal{H}\}$ as the active power at the sending end of branch connected to prosumer $i$. Similarly, the reactive power at the sending end of branch $i$ is denoted by $\bm{q}_i^{\text{injection}}= \{q_i^{\text{injection}}[t],~\forall t \in \mathcal{H}\}$. The reactive load at node $i$ is denoted as $\bm{q}_i^{\text{load}}= \{q_i^{\text{load}}[t],~\forall t \in \mathcal{H}\}$. Both active power and reactive power are upper-bounded and lower-bounded. We denote $\overline{P}_{\text{DN}}$ and $\underline{P}_{\text{DN}}$ as the upper and lower bounds for the active power at all times. Similarly, we denote $\overline{Q}_{\text{DN}}$ and $\underline{Q}_{\text{DN}}$ as the upper and lower bounds for the reactive power at all times. Hence, we have the following constraints for the active and reactive power:
\begin{align}
    & \rv{\bm{p}_{i+1}^{\text{injection}} = \bm{p}_i^{\text{injection}} + \bm{p}_i^{\text{output}},} \label{constraint-network1} \\
    & \rv{\overline{\bm{P}}_{\text{DN}} \leq \bm{p}_i^{\text{injection}} \leq \underline{\bm{P}}_{\text{DN}},} \label{constraint-network2}\\
    & \rv{\bm{q}_{i+1}^{\text{injection}} = \bm{q}_i^{\text{injection}} - \bm{q}_i^{\text{load}},} \label{constraint-network3} \\
    & \rv{\overline{\bm{Q}}_{\text{DN}} \leq \bm{q}_i^{\text{injection}} \leq \underline{\bm{Q}}_{\text{DN}},} \label{constraint-network4}
\end{align}
where $\overline{\bm{P}}_{\text{DN}}$, $\underline{\bm{P}}_{\text{DN}}$, $\overline{\bm{Q}}_{\text{DN}}$, and $\underline{\bm{Q}}_{\text{DN}}$ are the vector forms of $\overline{P}_{\text{DN}}$, $\underline{P}_{\text{DN}}$, $\overline{Q}_{\text{DN}}$ and $\underline{Q}_{\text{DN}}$, respectively.

The voltage of node $i$ is denoted as $\bm{V}_i = \{V_i [t],~\forall t \in \mathcal{H}\}$,
$r_i$ is the resistance of branch $i$, $x_i$ is the reactance of branch $i$, and $\epsilon$ is the tolerance parameter for voltage deviation. The voltage satisfies
\begin{align}
    & \rv{V_{i+1}[t] = V_i[t] - \frac{r_{i+1} p^{\text{injection}}_{i+1}[t] + x_{u+1} q^{\text{injection}}_{i+1}[t]} {V^{0}[t]},}~\forall t, \label{constraint-network5}\\
    & \rv{\bm{1}-\bm{\epsilon} \leq \bm{V}_{i} \leq \bm{1}+\bm{\epsilon},} \label{constraint-network6}
\end{align}
where $V^{0}[t]$ is the voltage at the root bus, $\epsilon$ is usually set to be $0.05$ per unit (p.u.), and $\bm{\epsilon}$ is the vector form of $\epsilon$.

\subsection{Optimization Problem Formulation}\label{sec:formulation}
Smart home prosumers cooperatively schedule their supply and demand and trade energy with each other to minimize the overall cost, consisting of grid energy purchase cost, battery degradation cost, thermal discomfort cost, energy trading cost, and revenue from feed-in solar energy.

We use a two-part electricity tariff for prosumer $i$'s grid energy purchase cost over the horizon $\mathcal{H}$ as
    \begin{align*}
        C_i^{\text{G}} = \pi_1^{\text{G}} \bm{1}^{\top} \bm{p}_i^{\text{G}} + \pi_2^{\text{G}} \norm{\bm{p}_i^{\text{G}}}_{\infty}, 
    \end{align*}
where $\top$ is transpose, $\pi_1^{\text{G}}$ is the energy rate, and $\pi_2^{\text{G}}$ is the peak demand rate. In addition to the traditional energy charge $\pi_1^{\text{G}} \bm{1}^{\top} \bm{p}_i^{\text{G}} = \pi_1^{\text{G}} \sum\nolimits_{t\in \mathcal{H}} p_i^{\text{G}}[t]$, this two-part electricity bill adds a peak charge $\pi_2^{\text{G}} \norm{\bm{p}_i^{\text{G}}}_{\infty} = \pi_2^{\text{G}} \max_{t\in \mathcal{H}} p_i^{\text{G}}[t]$ to shave the peak demand.

Prosumer $i$ can also feed in solar energy to the grid and gets paid at the feed-in tariff $\pi^{\text{feedin}}$. Hence, prosumer $i$'s revenue from selling solar energy is
    \begin{align*}
        R_i^{\text{S}} = \pi^{\text{feedin}} \bm{1}^{\top} \bm{p}_i^{\text{feedin}}.
    \end{align*}

The battery degradation cost of prosumer $i$ is modeled as a quadratic function of discharge energy as
    \begin{align*}
            C_i^{\text{B}}  = \beta^{\text{B}} \norm{\bm{p}_i^{\text{dis}}}_{2}^{2}, 
    \end{align*}
which shows that a deeper discharge will degrade the battery more, and $\beta^{\text{B}}$ is the degradation cost coefficient.

The thermal discomfort cost of smart home prosumer $i$ is measured by the deviation between the indoor temperature $\bm{T}_i^{\text{in}}$ and the set point $\bm{T}_i^{\text{REF}} = \{ T_i^{\text{REF}}[t],~\forall t \in \mathcal{H}\}$ as
    \begin{align*}
            C_i^{\text{hvac}}  = \beta^{\text{hvac}} \norm{\bm{T}_i^{\text{in}} - \bm{T}_i^{\text{ref}}}_{1}, 
    \end{align*}
where $\beta^{\text{hvac}}$ is the coefficient for the thermal discomfort cost.

We denote prosumer $i$'s energy scheduling as $\bm{p}_i^{\text{sch}}= \{ \bm{p}_i^{\text{G}}, \bm{p}_i^{\text{S}}, \bm{p}_i^{\text{feedin}}, \bm{p}_i^{\text{hvac}}, \bm{p}_i^{\text{ch}}, \bm{p}_i^{\text{dis}} \}$. We add up all the costs and revenues associated with energy scheduling decisions as the energy scheduling cost $C_i^{\text{sch}} (\bm{p}_i^{\text{sch}})$ as
    \begin{align*}
            C_i^{\text{sch}} (\bm{p}_i^{\text{sch}}) = C_i^{\text{G}} + C_i^{\text{B}} + C_i^{\text{hvac}} - R_i^{\text{S}}.
    \end{align*}

After formulating the energy scheduling cost, we model the energy trading cost $C_i^{\text{trade}}(\bm{p}_i^{\text{trade}})$. \rv{We assume that the prosumers are price takers, and the local network operator sets the energy trading prices $\bm{\pi}^{\text{trade}} = \{ \pi^{\text{trade}}[t],~\forall t \in \mathcal{H}\}$ as incentive signals to guide the energy trading, as the operator has a global view of the network. Since our work focuses on the algorithm design for energy trading, we take the energy trading prices as given.} The prosumer $i$'s energy trading cost is shown as
    \begin{align*}
            C_i^{\text{trade}} (\bm{p}_i^{\text{trade}}) 
         = \bm{\pi}^{\text{trade} \top}\bm{p}_i^{\text{trade}},
    \end{align*}
where the energy trading cost can be negative, meaning that the prosumer sells more and gains a profit.

Therefore, we formulate the cooperative energy management problem to minimize the overall cost of both energy scheduling cost $C_i^{\text{sch}}(\bm{p}_i^{\text{sch}})$ and energy trading cost $C_i^{\text{trade}}(\bm{p}_i^{\text{trade}})$ for all smart home prosumers.\\
\leftline{\textbf{Total Cost Minimization for Smart Home Prosumers}}
    \begin{align*}
        & \text{minimize} 
        && \sum\nolimits_{i\in\mathcal{N}} \left( C_i^{\text{sch}} (\bm{p}_i^{\text{sch}}) + C_i^{\text{trade}} (\bm{p}_i^{\text{trade}}) \right)  \\
        &\text{subject to} &&
        \text{Constraints}~\eqref{constraint-solar1}-\eqref{constraint-network6} \\
        &\text{Variables} &&
        \{ \bm{p}_i^{\text{sch}}, \bm{p}_{i}^{\text{trade}},\bm{p}_i^{\text{injection}},\bm{q}_i^{\text{injection}},\bm{V}_i,~\forall i \in \mathcal{N} \}.
    \end{align*}

\rv{Traditionally, the local network operator solves the total cost minimization problem, but some information of prosumers is private and may not be revealed to anyone else. Specially, internal operations, e.g., $\bm{p}_i^{\text{S}}$, $\bm{p}_i^{\text{hvac}}$, $\bm{p}_i^{\text{ch}}$, and $\bm{p}_i^{\text{dis}}$, do not go through the local network and is regarded as private information of smart home prosumers that should not be revealed to others. Other variables, such as grid purchase $\bm{p}_i^{\text{G}}$, feed-in solar energy $\bm{p}_i^{\text{feedin}}$, and energy trading $\bm{p}_i^{\text{trade}}$ are all visible to the operator. Therefore, we seek to solve the total cost minimization problem for smart home prosumers in a distributed manner to preserve their privacy and reduce the barrier to participation in energy trading. The distributed solution method is presented in Section~\ref{sec:solution}.}

\section{Distributed Solution Method}\label{sec:solution}
We develop an asynchronous distributed ADMM algorithm to solve the problem, suiting asynchronous communications and preserving the operational privacy of smart home prosumers. Designing an asynchronous distributed algorithm is motivated by the fact that different prosumers have different communication capacities and exchange information with the operator at different paces, thus demanding an asynchronous algorithm.

\subsection{\rv{Problem Decomposition}}\label{sec:decomposition}
We see from the total cost minimization problem for smart home prosumers that prosumers' energy trading decisions are coupled, and thus we need to introduce auxiliary variables $\tilde{\bm{p}}_{i}^{\text{trade}} = \{\tilde{p}_{i}^{\text{trade}}[t],~\forall t \in \mathcal{H}\}$, which should be identical to $\bm{p}_{i}^{\text{trade}}$. Hence, we can replace the energy trading balance constraint \eqref{constraint-trading1} by 
    \begin{align}
        & \tilde{\bm{p}}_{i}^{\text{trade}} = \bm{p}_{i}^{\text{trade}}, \label{constraint-auxiliary1} \\ %
        & \sum_{i \in \mathcal{N}} \tilde{\bm{p}}_{i}^{\text{trade}} = \bm{0}, \label{constraint-auxiliary2} 
    \end{align}
in which the auxiliary variables $\tilde{\bm{p}}_{i}^{\text{trade}}$ will be updated by the local network operator and provided to smart home prosumers as given parameters. Replacing the energy trading balance constraint \eqref{constraint-trading1} among prosumers by \eqref{constraint-auxiliary2} will help decompose the optimization problem to be solved by each prosumer in parallel. 

In addition to the energy trading $\bm{p}_{i}^{\text{trade}}$, prosumers' grid energy purchase and feed-in solar energy also go through the network and thus create coupling between prosumers and the network. Hence, we define another variable denoted as the net load $\bm{p}_i^{\text{net}} = \bm{p}_i^{\text{G}} - \bm{p}_i^{\text{feedin}}$ and introduce another auxiliary variable for the local network operator $\tilde{\bm{p}}_{i}^{\text{net}} = \{\tilde{p}_{i}^{\text{net}}[t],~\forall t \in \mathcal{H}\}$, which should be identical to $\bm{p}_{i}^{\text{net}}$. Similarly, we have \eqref{constraint-trading1} by 
    \begin{align}
        \tilde{\bm{p}}_{i}^{\text{net}} = \bm{p}_{i}^{\text{net}}, \label{constraint-auxiliary3}
    \end{align}
and replace the output power constraint in \eqref{constraint-output} as
\begin{align}
    \rv{\bm{p}_i^{\text{output}} = - \tilde{\bm{p}}_{i}^{\text{net}} - \tilde{\bm{p}}_i^{\text{trade}}.} ~\label{constraint-auxiliary4}
\end{align}

To decompose the coupling constraints of the auxiliary variables, we introduce dual variables $ \bm{\lambda}_i^{\text{trade}} = \{ \lambda_i^{\text{trade}} [t],~\forall t \in \mathcal{H}\}$ for Constraint \eqref{constraint-auxiliary1}, $\bm{\lambda}_i^{\text{net}} = \{ \lambda_i^{\text{net}} [t],~\forall t \in \mathcal{H}\}$ for Constraint \eqref{constraint-auxiliary3}, and obtain the augmented Lagrangian as
    \begin{align*}
        L &=  
        \sum_{i \in \mathcal{N}} 
        \Big(
        C_i (\bm{p}_i^{\text{sch}}) + C_i^{\text{trade}} (\bm{p}_i^{\text{trade}})  \\
        &+
        \bm{\lambda}_i^{\text{trade} \top} \left( \tilde{\bm{p}}_{i}^{\text{trade}} - \bm{p}_{i}^{\text{trade}} \right)
         + \frac{\rho_i^{\text{trade}}}{2} \left \| \tilde{\bm{p}}_{i}^{\text{trade}} - \bm{p}_{i}^{\text{trade}} \right \|^{2} \\
        &+
        \bm{\lambda}_i^{\text{net} \top} \left( \tilde{\bm{p}}_{i}^{\text{net}} - \bm{p}_{i}^{\text{net}} \right)
         + \frac{\rho_i^{\text{net}}}{2} \left \| \tilde{\bm{p}}_{i}^{\text{net}} - \bm{p}_{i}^{\text{net}} \right \|^{2} 
         \Big) ,
    \end{align*}
where $\rho_i^{\text{trade}} >0$ and $\rho_i^{\text{net}} >0$ are penalty parameters.
Based on the augmented Lagrangian, each smart home prosumer can solve the energy scheduling $\bm{p}_i^{\text{sch}}$ and energy trading $\bm{p}_{i}^{\text{trade}}$ as a local optimization problem based on the corresponding dual variables $\{ \bm{\lambda}_{i}^{\text{trade}}, \bm{\lambda}_{i}^{\text{net}}\}$ and auxiliary variables $\{ \tilde{\bm{p}}_{i}^{\text{trade}}, \tilde{\bm{p}}_{i}^{\text{net}} \}$, such that the operational privacy can be preserved.

\begin{figure}[t]
    \centering
    \includegraphics[width=0.9\linewidth]{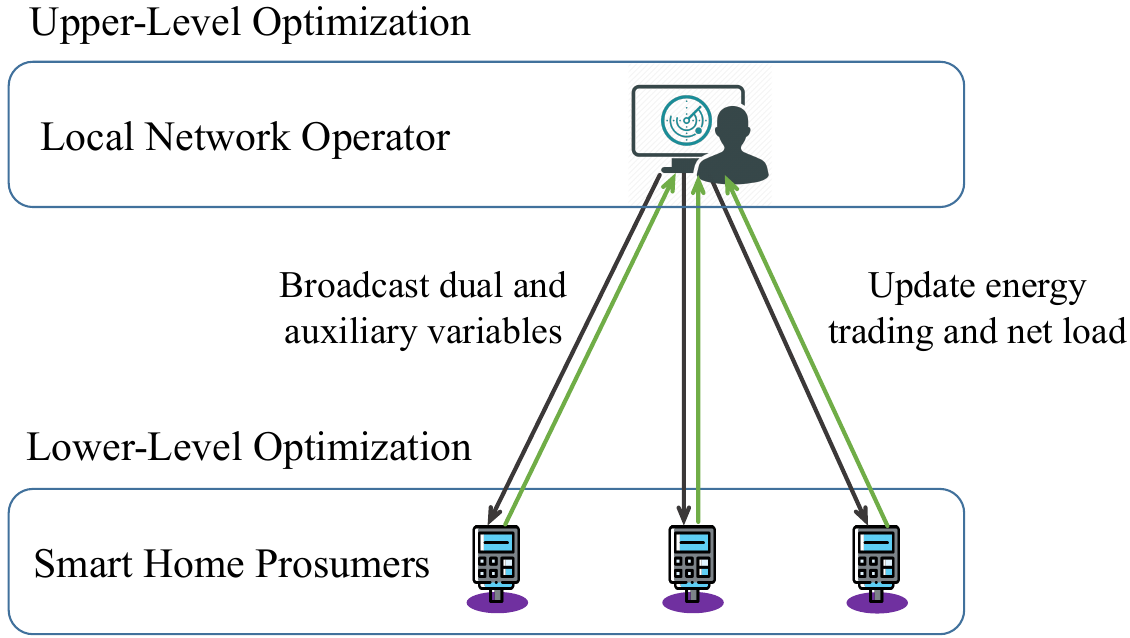}
    \caption{\rv{The upper-level and lower-level optimization problems for the interaction between the local network operator and smart home prosumers.}}
    \label{f:level}
\end{figure}

\rv{The total cost minimization for smart home prosumers can be decomposed into two subproblems, including an upper-level optimization problem and a lower-level optimization problem, as shown in Fig.~\ref{f:level}. The local network operator solves the dual variables $\{ \bm{\lambda}_{i}^{\text{trade}}, \bm{\lambda}_{i}^{\text{net}}\}$ and auxiliary variables $\{ \tilde{\bm{p}}_{i}^{\text{trade}}, \tilde{\bm{p}}_{i}^{\text{net}} \}$ based on the updates on energy trading and net load from prosumers and then broadcasts the dual variables and auxiliary variables to corresponding prosumers. To be specific, we first present the lower-level optimization problem for each prosumer as}

\leftline{\textbf{LP$_{i}$: Lower-Level Optimization for Prosumer $i$}}
    \begin{align*}
        &\text{minimize} && C_i^{\text{sch}} (\bm{p}_i^{\text{sch}}) + C_i^{\text{trade}} (\bm{p}_i^{\text{trade}}) \\
        &&& - \bm{\lambda}_i^{\text{trade} \top} \bm{p}_{i}^{\text{trade}}
         + \frac{\rho_i^{\text{trade}}}{2} \left \| \tilde{\bm{p}}_{i}^{\text{trade}} - \bm{p}_{i}^{\text{trade}} \right \|^{2}  \\
        &&& - \bm{\lambda}_i^{\text{net} \top} \bm{p}_{i}^{\text{net}}
         + \frac{\rho_i^{\text{net}}}{2} \left \| \tilde{\bm{p}}_{i}^{\text{net}} - \bm{p}_{i}^{\text{net}} \right \|^{2}  \\
        &\text{subject to} &&
        \text{Constraints}~\eqref{constraint-solar1}-\eqref{constraint-hvac2}~\text{and}~\eqref{constraint-trading2}-\eqref{constraint-balance}\\
        &\text{Variables} &&
        \{ \bm{p}_i^{\text{sch}}, \bm{p}_{i}^{\text{trade}}, \bm{p}_{i}^{\text{net}} \},
    \end{align*}
which solves the optimal energy scheduling including net load and energy trading by each prosumer $i$ independently in parallel. Note that the corresponding dual variables $\{ \bm{\lambda}_{i}^{\text{trade}}, \bm{\lambda}_{i}^{\text{net}} \}$ and auxiliary variables $\{ \tilde{\bm{p}}_{i}^{\text{trade}}, \tilde{\bm{p}}_{i}^{\text{net}} \}$ are given by the local network operator.

\rv{Then the local network operator solves the upper-level optimization problem as}

\leftline{\textbf{UP: Upper-Level Optimization for The Operator}}
    \begin{align*}
        &\text{minimize} && \sum_{i \in \mathcal{N}} 
        \Big( 
        \bm{\lambda}_i^{\text{trade} \top} \tilde{\bm{p}}_{i}^{\text{trade}}
         + \frac{\rho_i^{\text{trade}}}{2} \left \| \tilde{\bm{p}}_{i}^{\text{trade}} - \bm{p}_{i}^{\text{trade}} \right \|^{2} \\
        &&& + \bm{\lambda}_i^{\text{net} \top} \tilde{\bm{p}}_{i}^{\text{net}}
         + \frac{\rho_i^{\text{net}}}{2} \left \| \tilde{\bm{p}}_{i}^{\text{net}} - \bm{p}_{i}^{\text{net}} \right \|^{2} 
        \Big) \\
        &\text{subject to} &&
        \text{Constraints}~\eqref{constraint-network1}-\eqref{constraint-network6},~\eqref{constraint-auxiliary2}~\text{and}~\eqref{constraint-auxiliary4}\\
        &\text{Variables} &&
        \{ \tilde{\bm{p}}_{i}^{\text{trade}},\tilde{\bm{p}}_{i}^{\text{net}},\bm{p}_i^{\text{injection}},\bm{q}_i^{\text{injection}},\bm{V}_i,~\forall i \in \mathcal{N} \}.
    \end{align*}

The traditional ADMM algorithm solves the augmented Lagrangian in a synchronous manner, i.e., each iteration requires all smart home prosumers to update the energy trading decisions.
But it is difficult and inefficient to coordinate smart home prosumers using synchronous communications in practice, which motivates us to develop an asynchronous
distributed ADMM algorithm in Section~\ref{sec:algorithm}.

\subsection{Asynchronous
Distributed ADMM Algorithm Design}\label{sec:algorithm}
\begin{figure}[t]
    \centering
    \includegraphics[width=1.0\linewidth]{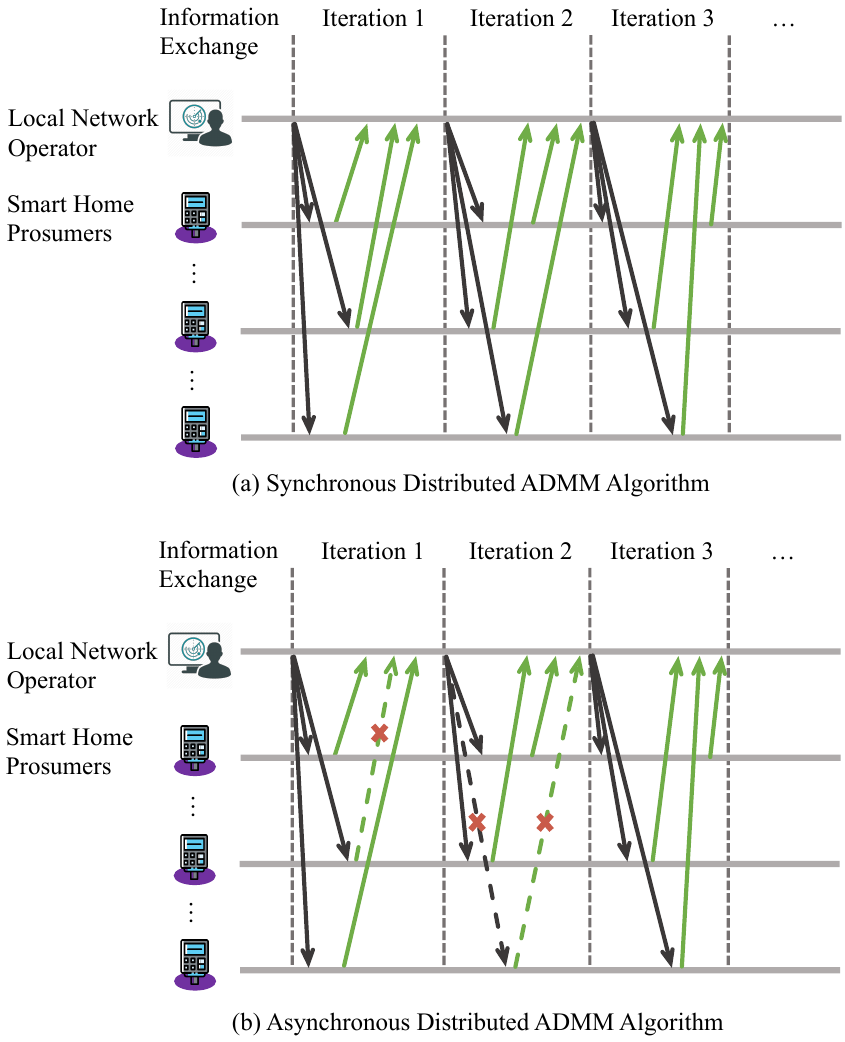}
    \caption{\rv{The illustration of information exchange in (a) the synchronous distributed ADMM algorithm and (b) the asynchronous distributed ADMM algorithm.}}
    \label{f:algorithm}
\end{figure}

Unlike traditional ADMM, we develop an asynchronous distributed ADMM algorithm \cite{kumar2016asynchronous} to solve the augmented Lagrangian, allowing asynchronous communications among smart home prosumers and thus working for heterogeneous communication capabilities. \rv{Fig.~\ref{f:algorithm} illustrates the difference between the synchronous distributed ADMM algorithm and the asynchronous distributed ADMM algorithm used in our work. As we see in Fig.~\ref{f:algorithm}(a), in each iteration, synchronous communication is needed to complete the information exchange between the local network operator and all prosumers. Then the operator will proceed to calculate the dual and auxiliary variables and make them ready to broadcast to prosumers in the next iteration. In contrast, we can see in Fig.~\ref{f:algorithm}(b) that, in some iterations, there could exist delayed or missing information exchange, shown as dash lines, due to prosumers' heterogeneous communication capabilities. The operator does not need to wait for the updates from all prosumers but proceeds to calculate dual and auxiliary variables. In the future iterations, those prosumers who fail to update in the previous iterations can update their energy trading net load decisions to the operator, which facilitates a more flexible information exchange paradigm.}

We denote $k$ as the iteration time and denote $\mathcal{N}(k)$ as the set of smart home prosumers that have received the corresponding dual variables $\{ \bm{\lambda}_{i}^{\text{trade}}(k), \bm{\lambda}_{i}^{\text{net}}(k) \}$ and auxiliary variables $\{ \tilde{\bm{p}}_{i}^{\text{trade}}(k), \tilde{\bm{p}}_{i}^{\text{net}}(k) \}$ and updated their energy trading and net load decisions, which are ready to send to the operator. For other prosumers in $ \mathcal{N} \backslash \mathcal{N}(k)$, they either do not receive the corresponding dual and auxiliary variables or do not update new decisions in this iteration. In this case, the local network operator does not need to wait for all prosumers' updates but uses the previous updated decisions for those who do not update in this iteration. Also, those prosumers who do not update can wait for the receipt of dual variables and auxiliary variables in future iterations to update their decisions to the operator.
Therefore, in iteration $k$, smart home prosumer $i$ solves energy scheduling and trading decisions as $\bm{p}_i^{\text{sch}}(k+1), \bm{p}_{i}^{\text{trade}}(k+1)$. Only the energy trading and net load decisions $\{ \bm{p}_{i}^{\text{trade}}(k+1), \bm{p}_{i}^{\text{net}}(k+1) \}$ will be sent to the operator. For simplicity, we denote the objective function of the Lower-Level Optimization for Prosumer $i$, known as \textbf{LP$_{i}$}, as 
    \begin{align*}
        C_i^{\text{LP}} & = C_i^{\text{sch}} (\bm{p}_i^{\text{sch}}) + C_i^{\text{trade}} (\bm{p}_i^{\text{trade}}) \\
        & - \bm{\lambda}_i^{\text{trade} \top} \bm{p}_{i}^{\text{trade}}
         + \frac{\rho_i^{\text{trade}}}{2} \left \| \tilde{\bm{p}}_{i}^{\text{trade}} - \bm{p}_{i}^{\text{trade}} \right \|^{2}  \\
        & - \bm{\lambda}_i^{\text{net} \top} \bm{p}_{i}^{\text{net}}
         + \frac{\rho_i^{\text{net}}}{2} \left \| \tilde{\bm{p}}_{i}^{\text{net}} - \bm{p}_{i}^{\text{net}} \right \|^{2}.
    \end{align*}

Due to asynchronous communications, in iteration $k+1$, some prosumers' updates are received by the operator, but other prosumers in $\mathcal{N} \backslash \mathcal{N}(k)$ do not successfully update their decisions. Hence, the operator takes the following values for prosumers' energy trading and net load as
    \begin{equation}
    \begin{aligned}
        &\{ \bm{p}_i^{\text{trade}}(k+1), \bm{p}_{i}^{\text{net}}(k+1) \} = \\
        &\left\{
          \begin{array}{rcl}
          \arg \min C_i^{\text{LP}}, && {i \in \mathcal{N}^k},\\
          \{ \bm{p}_i^{\text{trade}}(k), \bm{p}_{i}^{\text{net}}(k) \} && {i \in \mathcal{N} \backslash \mathcal{N}^k}.
          \end{array} 
          \right. \label{update-lp}
    \end{aligned}
    \end{equation}
The rest of the prosumers' decisions, such as battery scheduling and HVAC control, are regarded as operational privacy for prosumers, which will not be revealed to anyone else, even the operator.

The local network operator, in the upper-level optimization problem \textbf{UP}, solves the auxiliary variables and also updates the dual variables for all smart home prosumers. Specifically, the operator solves the \textbf{UP} to obtain the auxiliary variables as
    \begin{equation}
    \begin{aligned}
        & \{ \tilde{\bm{p}}_{i}^{\text{trade}}(k), \tilde{\bm{p}}_{i}^{\text{net}}(k) \} = \\
        & \arg \min
        \sum_{i \in \mathcal{N}} 
        \Big( 
        \bm{\lambda}_i^{\text{trade} \top} \tilde{\bm{p}}_{i}^{\text{trade}}
         + \frac{\rho_i^{\text{trade}}}{2} \left \| \tilde{\bm{p}}_{i}^{\text{trade}} - \bm{p}_{i}^{\text{trade}} \right \|^{2} \\
        & + \bm{\lambda}_i^{\text{net} \top} \tilde{\bm{p}}_{i}^{\text{net}}
         + \frac{\rho_i^{\text{net}}}{2} \left \| \tilde{\bm{p}}_{i}^{\text{net}} - \bm{p}_{i}^{\text{net}} \right \|^{2} 
        \Big), \label{update-up}
    \end{aligned}
    \end{equation}
subject to the network constraints \eqref{constraint-network1}-\eqref{constraint-network6} and constraints for auxiliary variables, including  \eqref{constraint-auxiliary2} and \eqref{constraint-auxiliary4}.

After solving the auxiliary variables $\{ \tilde{\bm{p}}_{i}^{\text{trade}}(k), \tilde{\bm{p}}_{i}^{\text{net}}(k) \}$, the local network operator further updates the dual variables as
     \begin{align}
        & \bm{\lambda}_{i}^{\text{trade}}(k+1) \leftarrow \bm{\lambda}_{i}^{\text{trade}}(k) + \rho_{i}^{\text{trade}} \left( \tilde{\bm{p}}_{i}^{\text{trade}}(k) - \bm{p}_{i}^{\text{trade}}(k) \right), \label{update-dual_trade}  \\
        & \bm{\lambda}_{i}^{\text{net}}(k+1) \leftarrow \bm{\lambda}_{i}^{\text{net}}(k) + \rho_{i}^{\text{net}} \left( \tilde{\bm{p}}_{i}^{\text{net}}(k) - \bm{p}_{i}^{\text{net}}(k) \right). \label{update-dual_net}
     \end{align}

Both updated auxiliary variables $\{ \tilde{\bm{p}}_{i}^{\text{trade}}(k+1), \tilde{\bm{p}}_{i}^{\text{net}}(k+1) \}$ and dual variables $\{ \bm{\lambda}_{i}^{\text{trade}}(k+1), \bm{\lambda}_{i}^{\text{net}}(k+1) \}$
which will be sent to smart home prosumers in LP$_{i}$ to update their energy scheduling and trading in the next iteration. The entire process of our developed Asynchronous Distributed ADMM Algorithm is presented in Algorithm~\ref{alg1}.

\begin{algorithm}[t]
     \caption{\rv{Asynchronous Distributed ADMM Algorithm for Energy Trading}}
     \label{alg1} 
     \SetAlgoLined
     \textbf{Initialization}:\\
     iteration index $k {\leftarrow} 1$;\\
     convergence threshold $\epsilon_{1}, \epsilon_{2} {\leftarrow} 1.0 \times 10^{-2}$; \\
     iteration stepsize $\rho_i^{\text{trade}}(0),\rho_i^{\text{net}}(0) {\leftarrow} 1$; \\
     dual variables $\{ \bm{\lambda}_{i}^{\text{trade}}, \bm{\lambda}_{i}^{\text{net}} \} {\leftarrow} \bm{0}$;
    
    \While{$\sum_{i \in \bm{N}} \parallel \{ \bm{p}_i^{\text{trade}}(k), \bm{p}_{i}^{\text{net}}(k) \} - \{ \tilde{\bm{p}}_i^{\text{trade}}(k), \tilde{\bm{p}_{i}^{\text{net}}}(k) \} \parallel > \epsilon_{1} $ \\
    $\lor$ $\parallel \{ \Delta \bm{\lambda}^{\text{trade}}(k), \Delta \bm{\lambda}^{\text{net}}(k) \} \parallel > \epsilon_{2}$}{
        \For{$i \in \mathcal{N}$}{
        (1) Prosumers who received auxiliary and dual variables solve the Lower-Level Optimization problem and update energy trading and net load decisions $\{ \bm{p}_{i}^{\text{trade}}(k), \bm{p}_{i}^{\text{net}}(k) \}$;
        }
        
    (2) The local network operator uses prosumers' updates according to \eqref{update-lp};
     
    (3) The local network operator in \textbf{UP} updates the auxiliary variables $\{\tilde{\bm{p}}_{i}^{\text{trade}},\tilde{\bm{p}}_{i}^{\text{net}} \}$ according to \eqref{update-up};
    
    (4) The local network operator updates the dual variable $\{ \bm{\lambda}_{i}^{\text{trade}}, \bm{\lambda}_{i}^{\text{net}} \}$ are updated according to \eqref{update-dual_trade} and \eqref{update-dual_net};
    
    (5) The local network operator broadcasts both auxiliary variables and dual variables to the corresponding prosumers $i$;
    
    (6) Update the iteration time $k \leftarrow k+1$;
    }
 \textbf{Outputs}:
 Optimal energy scheduling and trading solutions $\{ \bm{p}_i^{\text{sch} \star}, \bm{p}_i^{\text{trade} \star},~\forall i \in \mathcal{N} \}$.
\end{algorithm}

Our developed asynchronous distributed ADMM algorithm solves prosumers' energy scheduling and energy trading in an iterative fashion.
The terminating condition depends on the residuals of primal variables and dual variables. Specifically, two thresholds, e.g., $\epsilon_{1}$ and $\epsilon_{2}$, are introduced for the residuals of primal variables (i.e., $\{ \bm{p}_i^{\text{trade}}(k), \bm{p}_{i}^{\text{net}}(k) \}$) comparing to auxiliary variables (i.e., $\{ \tilde{\bm{p}}_i^{\text{trade}}(k), \tilde{\bm{p}_{i}^{\text{net}}}(k) \}$) and for residuals of dual variables $\{ \Delta \bm{\lambda}^{\text{trade}}(k), \Delta \bm{\lambda}^{\text{net}}(k) \}$, where $\Delta \bm{\lambda}^{\text{trade}}(k) = \bm{\lambda}^{\text{trade}}(k) - \bm{\lambda}^{\text{trade}}(k-1)$ and $\Delta \bm{\lambda}^{\text{net}}(k) = \bm{\lambda}^{\text{net}}(k) - \bm{\lambda}^{\text{net}}(k-1)$.
Despite asynchronous updates of $\{ \bm{p}_i^{\text{trade}}, \bm{p}_{i}^{\text{net}} \}$,
the algorithm converges to the optimal solution of the original total cost minimization problem, if the stepsizes are set as $\rho_{i}^{\text{trade}}(k) = \rho_{i}^{\text{net}}(k) = 1/k$ \cite{kumar2016asynchronous}.

\section{Performance Evaluation}\label{sec:eval}
We validate our proposed asynchronous distributed ADMM algorithm for smart home prosumers behind a common low-voltage transformer using parameter configurations in \cite{yang2021privacy} and real-world solar generation, and load profiles \cite{wang2016cooperative,pecan}. 
To simulate the asynchronous communications, in each iteration, we randomly select a group (e.g., $20\%$) of prosumers that do not receive the auxiliary variables and dual variables or do not update their energy scheduling and trading decisions in time. 
We choose the number of prosumers as $10$ and $50$ to show the convergence of the algorithm. Some of the key parameters in the simulations are summarized as follows. For HVAC, we set $C_i = 3.3$ and $R_i = 1.35$. The upper and lower bounds for the indoor temperature are $15$ and $32$ degree Celsius, respectively. The coefficient for the thermal discomfort cost is $0.25$. The coefficients of charging and discharging efficiency are $0.9$. The battery degradation cost coefficient is set as $0.01$. For the two-part electricity tariff, the energy rate and peak demand rate are $0.2$ and $1.2$, respectively. 

\subsection{\rv{Convergence of The Asynchronous Algorithm}}\label{sec:convergence}
\begin{figure}[t] 
  \centering 
  \includegraphics[width=1.0\linewidth]{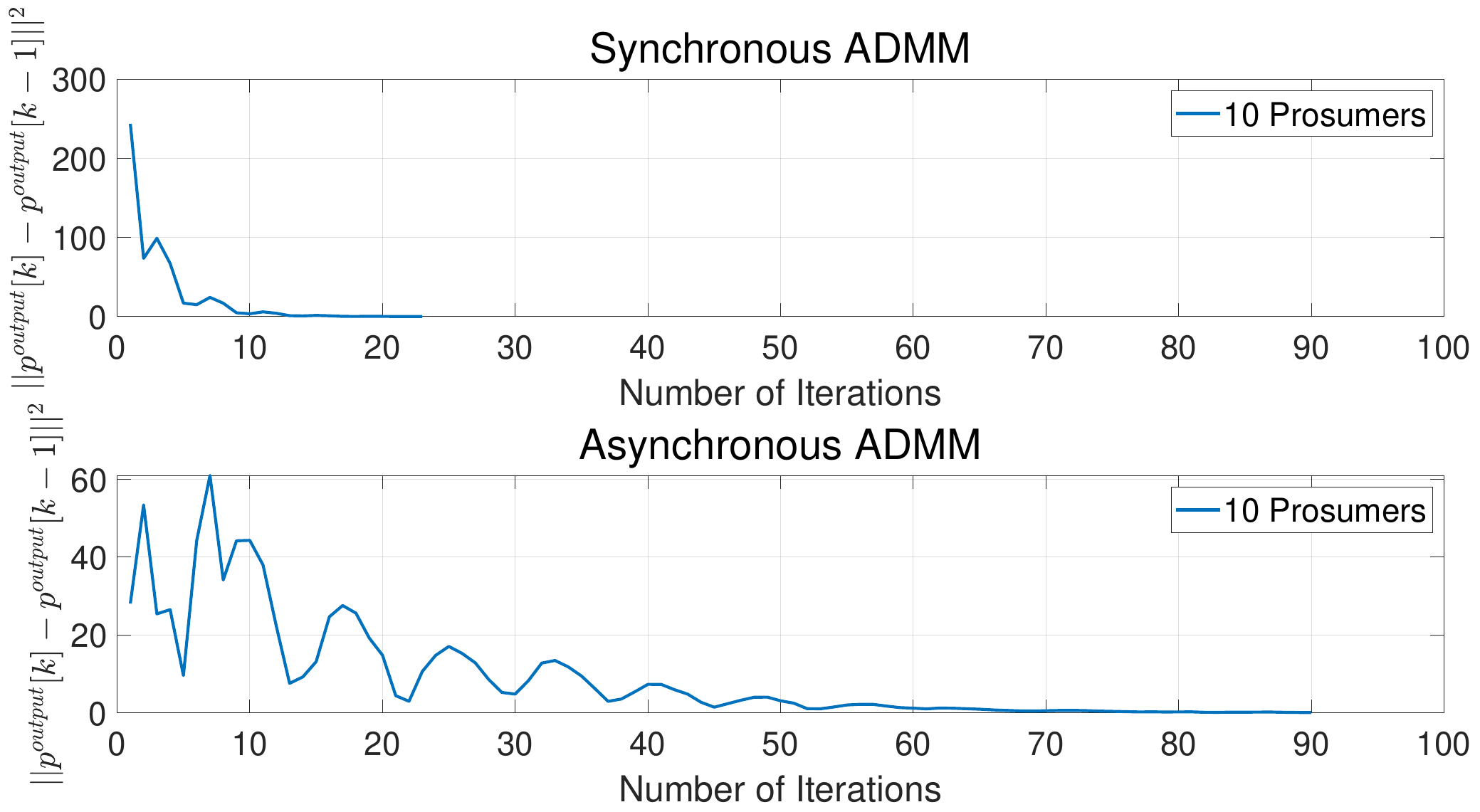}
  \caption{Comparison of the convergence of synchronous and asynchronous distributed ADMM algorithms with $10$ prosumers.}
  \label{f:convergence10} 
\end{figure}

We first validate the convergence of our developed asynchronous distributed ADMM algorithm with $10$ and $50$ prosumers. Fig.~\ref{f:convergence10} shows the comparison of the convergence of the standard synchronous ADMM algorithm and our asynchronous ADMM algorithm with $10$ prosumers. 
The X-axis indicates the number of iterations, and the Y-axis indicates convergence errors.
We see that the synchronous algorithm converges to the terminating threshold $1.0\times10^{-1}$ within $23$ iterations. In contrast, the asynchronous algorithm converges within $90$ iterations. Note that it is understandable that the asynchronous algorithm takes more iterations to converge. But in terms of time, the asynchronous algorithm could outperform the synchronous algorithm, as the asynchronous algorithm does not need to wait for all the updates from prosumers to complete each iteration\footnote{We need to set up a parallel computing environment for testing the running time of the two algorithms, which is beyond the scope of this work and considered as future work.}.

\begin{figure}[t] 
  \centering 
  \includegraphics[width=1.0\linewidth]{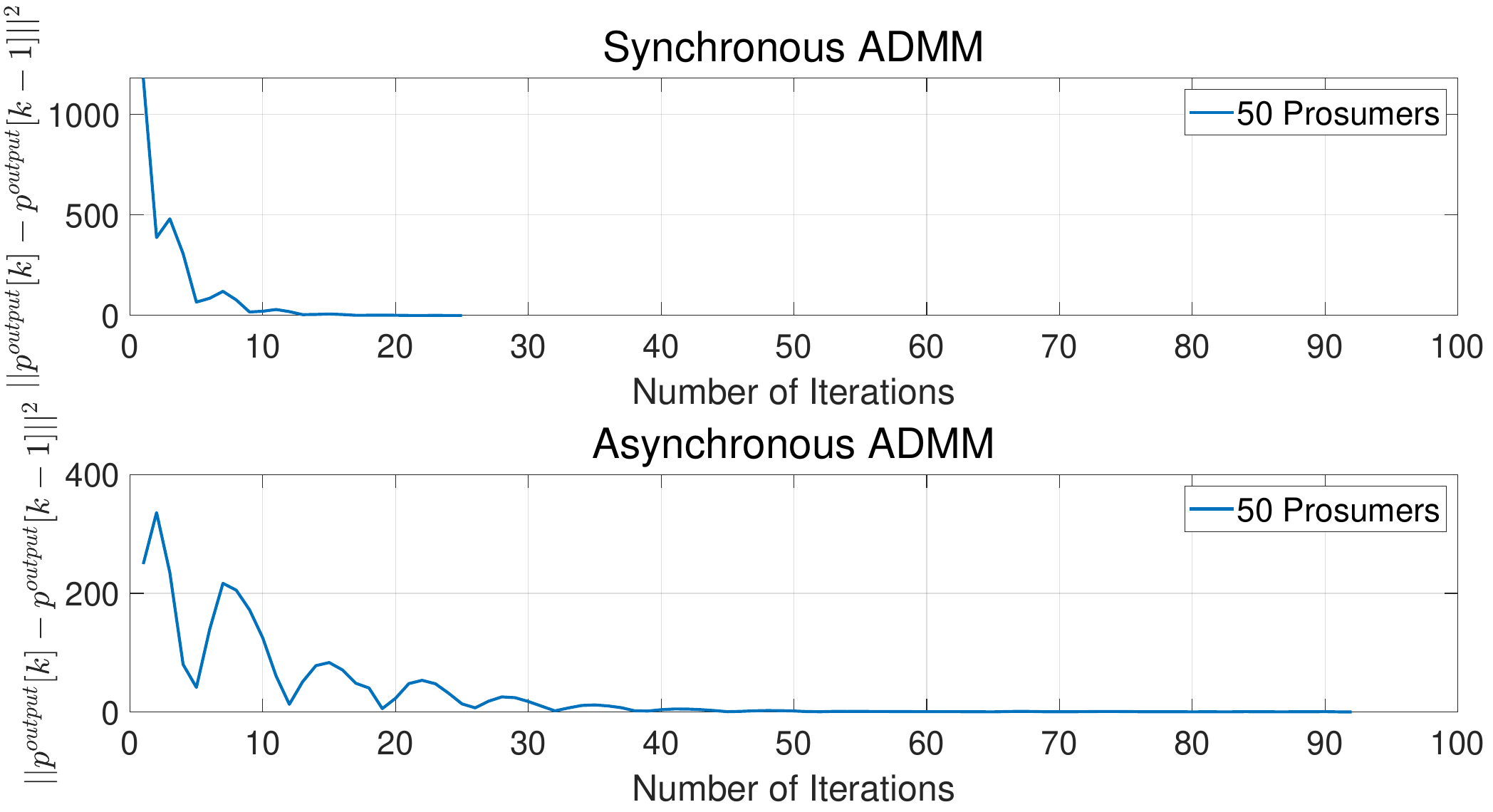}
  \caption{Comparison of the convergence of synchronous and asynchronous distributed ADMM algorithms with $50$ prosumers.}
  \label{f:convergence50} 
\end{figure}

We increase the number of prosumers to $50$, and Fig.~\ref{f:convergence50} shows the comparison of the convergence of the standard synchronous ADMM algorithm and our asynchronous ADMM algorithm with $50$ prosumers. We see that both algorithms can scale up. The synchronous algorithm converges within $25$ iterations, and the asynchronous algorithm converges within $92$ iterations.

\subsection{\rv{Optimal Energy Trading and Cost Reduction}}\label{sec:trade}
\begin{figure}[t]
   \centering 
    \includegraphics[width=1.0\linewidth]{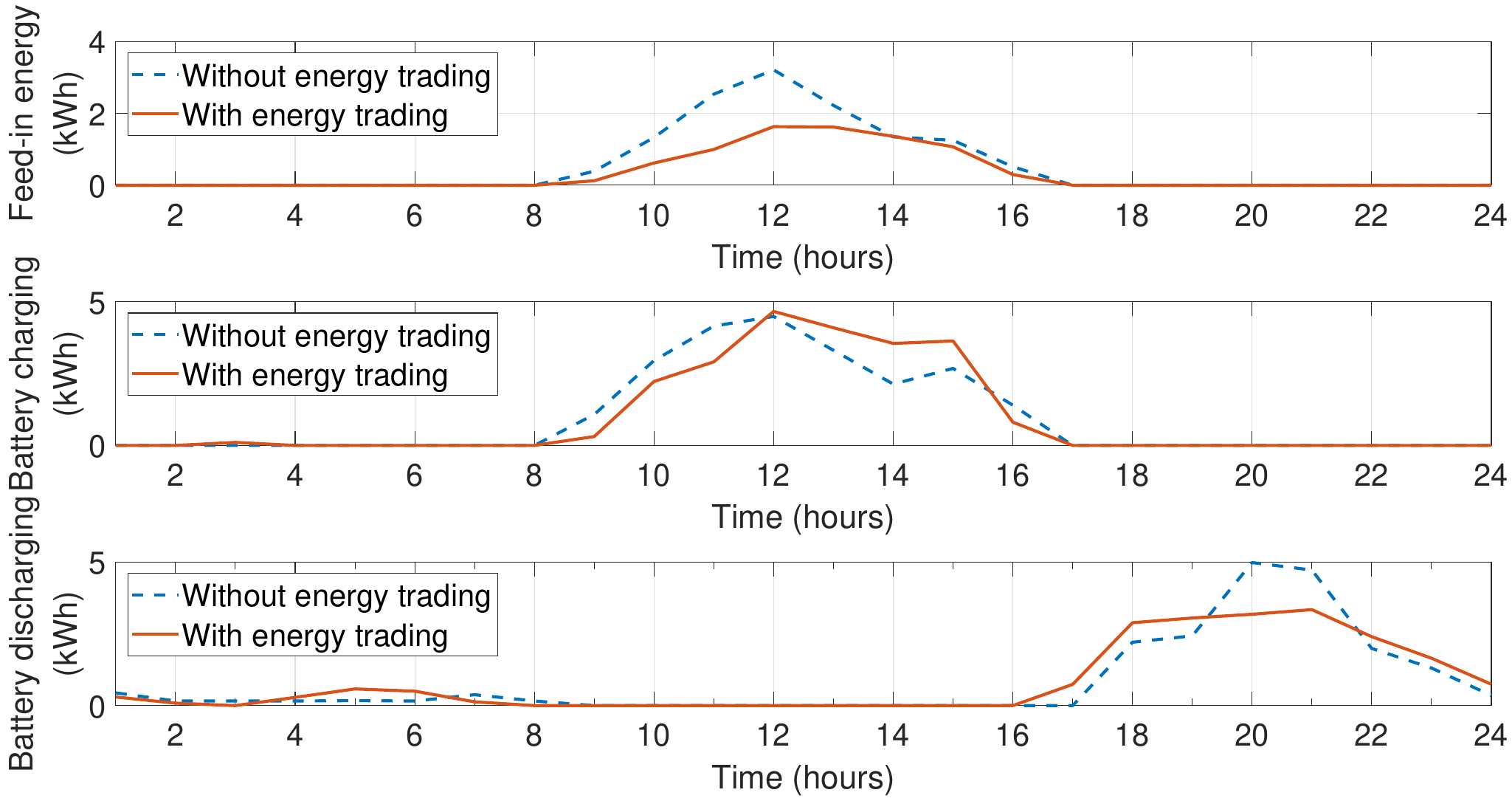} 
  \caption{Comparison of the optimal energy schedule of a representative prosumer, including feed-in solar energy, battery charge and discharge, without and with energy trading.} 
  \label{f:schedule} 
\end{figure}

\begin{figure}[t] 
   \centering 
    \includegraphics[width=1.0\linewidth]{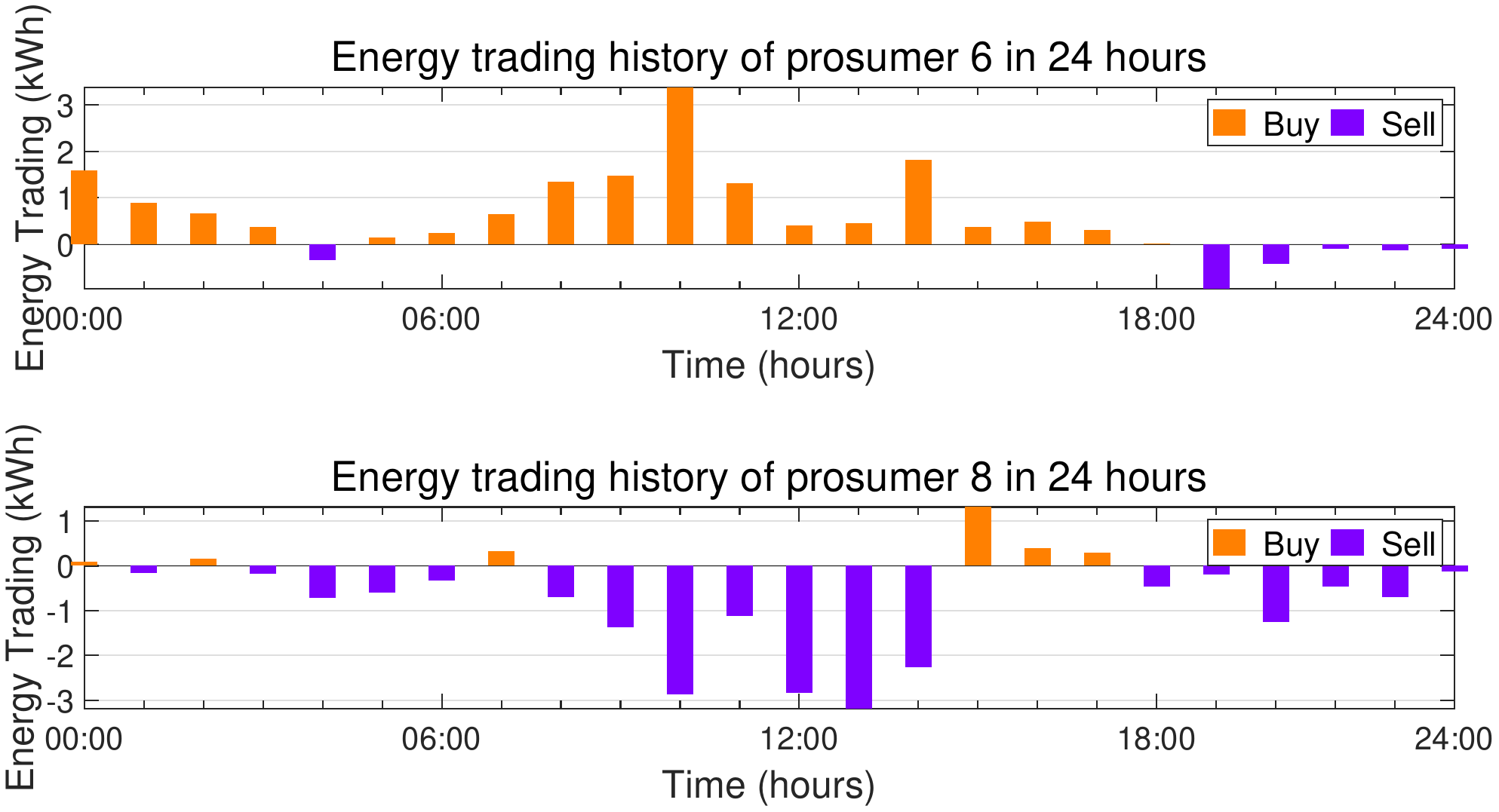} 
  \caption{Comparison of the optimal energy trading of two representative prosumers (6 and 8).} 
  \label{f:trade} 
\end{figure}

\begin{figure}[t] 
   \centering 
    \includegraphics[width=1.0\linewidth]{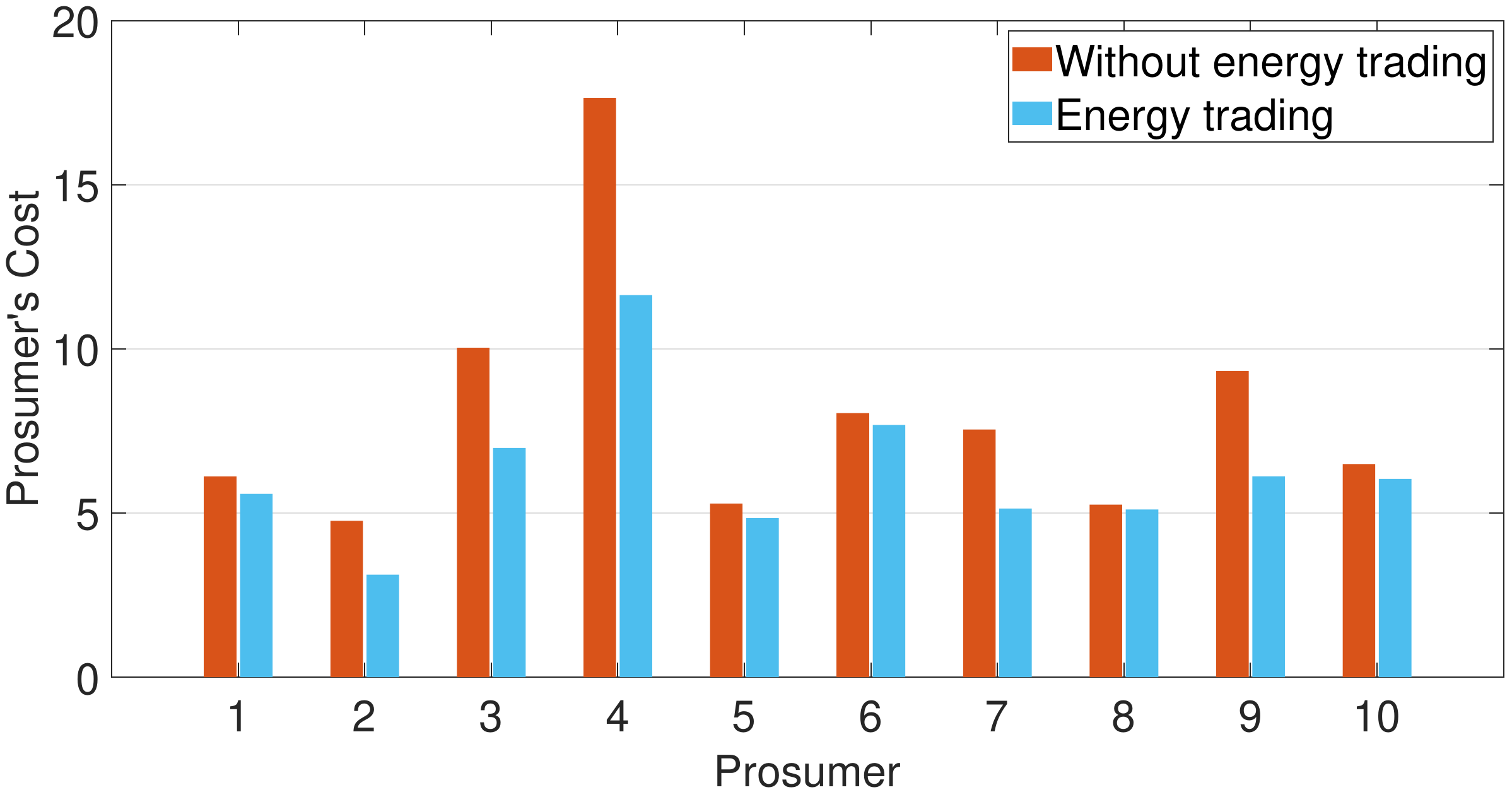} 
  \caption{Comparison of the minimized energy cost of all prosumers without and with energy trading.} 
  \label{f:costreduction} 
\end{figure}

We validate the effectiveness of energy trading for smart home prosumers. 
Fig.~\ref{f:schedule} depicts the comparison of the optimal energy schedule of a representative prosumer. We see that, without energy trading, the prosumer feeds in more solar energy to the grid. But with energy trading, the prosumer lowers the feed-in solar energy but charge more into the battery to discharge for energy trading later. The optimal battery charge and discharge profiles are slightly different due to the optimal energy trading among all prosumers in the local network.
Fig.~\ref{f:trade} shows the comparison of the optimal energy trading of two representative prosumers (6 and 8), exhibiting diverse trading profiles in different time slots. The orange bar indicates buying energy, and the purple bar indicates selling energy. Prosumer $6$ buys more, while Prosumer $8$ sells more in the local network, achieving a win-win situation to participating prosumers for mutual benefits.
Fig.~\ref{f:costreduction} shows the comparison of the minimized cost of all $10$ prosumers without and with energy trading. The red bar shows the minimized cost of individual prosumer without energy trading, and the blue bar shows the minimized cost of individual prosumer with energy trading. We see that individual prosumers' costs are reduced by up to $34\%$. The overall cost of all prosumers is reduced by $11.3\%$.
These results demonstrate that the proposed energy trading algorithm effectively leverages the diversity of prosumers to trade energy and reduce costs.

\section{Conclusion and Future Work}\label{sec:conclusion}
This paper developed a network-aware asynchronous distributed ADMM algorithm for energy management and energy trading for smart home prosumers, providing a more practical solution toward real-world application under asynchronous communications. Our developed algorithm iterated between an upper-level problem solved by the local network operator to meet the network constraints and many lower-level problems solved by prosumers in parallel. As two-way communications play a crucial role in managing DERs in the future power grid, our developed asynchronous algorithm can benefit the development of distributed coordination for heterogeneous DERs under realistic communication conditions.

For future work, we identify several critical issues that are worth investigating. First, we will develop more efficient 
asynchronous algorithms for energy trading that can scale up to a more practical scenario with more prosumers. Second, we will study the energy trading pricing mechanism reflecting the need and value of energy trading in real-time. Third, we will consider a larger distribution network to host more prosumers for simulations.

\bibliographystyle{IEEEtran}
\bibliography{IEEEabrv.bib,ref.bib}

\end{document}